\documentclass{article}
\begin{document}
\title{Well-posedness of boundary layer equations for time-dependent flow of non-Newtonian fluids}
\author{Michael Renardy \and Xiaojun Wang}

\maketitle

\begin{abstract}
We consider the flow of an upper convected Maxwell fluid in the limit of high Weissenberg and Reynolds number. In this limit, the no-slip
condition cannot be imposed on the solutions. We derive equations for the resulting boundary layer and prove the well-posedness of these
equations. A transformation to Lagrangian coordinates is crucial in the argument.
\end{abstract}

\def\err{{\rm I\!R}}
\def\T{{\bf T}}
\def\A{{\bf A}}
\def\B{{\bf B}}
\def\C{{\bf C}}
\def\F{{\bf F}}
\def\I{{\bf I}}
\def\Q{{\bf Q}}
\def\W{{\bf Wi}}
\def\R{{\bf Re}}
\def\v{{\bf v}}
\def\divg{\,{\rm div}\,}
\def\detg{\,{\rm det}\,}
\def\trg{\,{\rm tr}\,}

\def\y{{\bf y}}
\def\x{{\bf x}}
\def\n{{\bf n}}
\def\f{{\bf f}}
\def\0{{\bf 0}}
\def\S{{\bf S}}

\section{Introduction}

Classical fluid mechanics is based on the Navier-Stokes equations supplemented by a no-slip boundary condition on walls. In the limit
of zero viscosity, the Euler equations are obtained. However, the Euler equations do not allow for a no-slip boundary condition, and only
the normal component of the velocity can be prescribed to be zero on a wall. It was Prandtl's fundamental insight more than a century ago
\cite{PL} that, for many high Reynolds number flows, the Euler equations provide an adequate description except in a thin layer
close to the boundary, which is called a boundary layer. By taking advantage of the thinness of this layer, the Navier-Stokes equations
can formally be reduced to the system which is now known as the Prandtl equations. For two-dimensional flow, and a boundary placed
at $y=0$, these equations take the form
\begin{eqnarray}\label{prandtl}
{\partial u\over \partial t}+u{\partial u\over \partial x}+v{\partial u\over \partial y}&=& \nu u_{yy}-{\partial p\over \partial x},\nonumber \\
{\partial p\over \partial y}&=&0, \nonumber\\
{\partial u\over \partial x}+{\partial v\over \partial y}&=&0,\nonumber\\
u(x,y,0)&=&u_0(x,y), \nonumber\\
u(x,0,t)=v(x,0,t)&=&0,\nonumber\\
u(x,\infty,t)&=&u^{\infty}(x,0,t).
\end{eqnarray}
Here $u,v$ are velocities in $x,y$ directions, $p$ is pressure, $\nu$ is viscosity, and $u^{\infty}(x,y,t)$ represents the given flow in the core region.

One might hope to obtain a simplified procedure for solving high Reynolds number flow problems by solving the Euler equations (or even the simpler
special case of potential flow) in the core of the flow domain, and then solving the Prandtl equations near the boundary. This
program, however, runs into difficulties related to the question of well-posedness of the Prandtl equations. Oleinik \cite{Oleinik} established a well-posedness result under the assumption that the velocity profile in the boundary layer is monotone.
Sammartino and Caflisch \cite{Samm} established an existence result for analytic initial data. We also refer to the review article of Weinan E \cite{E} for further work prior to 2000. Recently, G\'erard-Varet and Dormy \cite{GV} established that, for general initial data, the
Prandtl equations are not well-posed in Sobolev spaces.

Viscoelastic flows exhibit phenomena of instability which share many characteristics of turbulence \cite{GS}. In this paper, we shall focus on the upper convected Maxwell model for the viscoelastic flow. In the limit
of high elasticity, a limiting equation can be derived which is similar to the system of ideal magnetohydrodynamics and, like the Euler equations, does not allow the imposition of a no-slip boundary condition. The well-posedness of this system has been
established in \cite{Wang}. The goal of this manuscript is to supplement this analysis with a study of the well-posedness of the
accompanying boundary layer equations.

We note that the ill-posedness of the Prandtl equations \cite{GV} is linked to shear flow instabilities. Elasticity has a
stabilizing effect on high Reynolds number flow instabilities \cite{hughes,ogilvie,miller,kaf,mr1}. This stabilizing
effect can restore well-posedness of the hydrostatic approximation in situations where this approximation is ill-posed
for the Euler equations \cite{mr2}. We may therefore hope to establish well-posedness in the boundary layer system as well.
Indeed, this turns out to be the case. We shall show that a transformation into Lagrangian coordinates transforms the
boundary layer system into a semilinear wave equation for which well-posedness can be readily established.

\section{The high Weissenberg number limit and the issue of boundary layers}

We start with the upper convected Maxwell model in dimensionless form:
\begin{eqnarray}\label{ucm1}
{\partial \v\over \partial t}+(\v\cdot\nabla) \v&=&{1\over R}\nabla\cdot \T-\nabla p, \nonumber \\
{\partial \T\over \partial t}+(\v\cdot\nabla) \T-(\nabla \v)\T-\T(\nabla \v)^T+{1\over W}\T&=&{1\over W}(\nabla \v+(\nabla \v)^T).
\end{eqnarray}
where $R$ is the Reynolds number and $W$ is the Weissenberg number measuring the elasticity.
Since normal stresses in shear flow are of order $W$ rather than order 1, we
shall also scale the stresses with an additional factor $W$ and obtain
\begin{eqnarray}\label{constitutive0}
{\partial \v\over \partial t}+(\v\cdot\nabla) \v&=&E\nabla\cdot \T-\nabla p,\nonumber\\
{\partial \T\over \partial t}+(\v\cdot\nabla) \T-(\nabla \v)\T-\T(\nabla \v)^T+{1\over W}\T&=&{1\over W^2}(\nabla \v+(\nabla \v)^T).
\end{eqnarray}
Here $E=W/R$ is the elasticity number. We are interested in a limit where $W$ and $R$ tend to infinity simultaneously, but $E$
remains fixed.

If we formally set $W=\infty$ above, we obtain the limiting system
\begin{eqnarray}\label{main1}
{\partial \v^0\over \partial t}+(\v^0\cdot\nabla) \v^0=E\nabla\cdot \T^0-\nabla p^0,\nonumber\\
{\partial \T^0\over \partial t}+(\v^0\cdot\nabla) \T^0-(\nabla \v^0)\T^0-\T^0(\nabla \v^0)^T={\bf 0},\nonumber\\
\nabla\cdot\v^0=0.
\end{eqnarray}
In \cite{Wang}, we proved the well-posedness of this system. We considered the initial-boundary value problem in a smooth domain
$\Omega$, subject to initial conditions for $\v^0$ and $\T^0$, and the boundary condition $\v^0\cdot\n=0$. A crucial assumption
was that $\T^0\cdot\n={\bf 0}$; it can be shown that the equations preserve this condition if it is satisfied initially. The physical background behind this assumption is that the local flow near a solid wall is always a shear flow, and at high Weissenberg number the extra stress is dominated by the first normal stress, i.e. the stress component tangent to the wall.

For the full equations, however, we have the boundary condition $\v={\bf 0}$, not just $\v\cdot\n=0$. To accommodate this, boundary layers must form near the wall. For discussing these boundary layers, it is convenient to set $\S=\T+\I/W^2$. With this substitution, the constitutive law transforms to
\begin{eqnarray}\label{constitutive}
{\partial \S\over \partial t}+(\v\cdot\nabla) \S-(\nabla \v) \S-\S(\nabla \v)^T+{1\over W}(\S-{1\over W^2}\I)&=&0,\nonumber
\end{eqnarray}
and we have $\divg\T=\divg\S$ in the momentum equation.
The change in boundary conditions is related to the fact that, while $\S\cdot\n$ vanishes on the boundary at leading order,
for the full equations we have strict positive definiteness of $\S$.

If we assume a flat boundary given by $y=0$, we can obtain a consistent scaling if we set
\begin{eqnarray}
t'&=&t,x'=x,y'=Wy,u'=u,v'=Wv,\nonumber\\
S_{11}'&=&S_{11},S_{12}'=WS_{12},S_{22}'=W^2S_{22},p'=p.
\end{eqnarray}
If we keep only the leading order terms and omit the primes in the rescaled equations, then the momentum equations reduce to
\begin{eqnarray}\label{scaling1}
{\partial u\over \partial t}+u{\partial u\over \partial x}+v{\partial u\over \partial y}=E({\partial S_{11}\over \partial x}+{\partial S_{12}\over \partial y})-{\partial p\over\partial x},\nonumber\\
{\partial p\over\partial y}=0.
\end{eqnarray}
The equation $\partial p/\partial y=0$ implies that  to the leading order the pressure $p(x,t)=p^0(x,0,t)$ in the boundary layer is the known function given by the outside flow.
For the constitutive equation and continuity equations, we have
\begin{eqnarray}\label{scaling2}
{\partial \S\over \partial t}+(\v\cdot\nabla) \S-(\nabla \v)\S-\S(\nabla \v)^T=0,\nonumber\\
{\partial u\over \partial x}+{\partial v\over \partial y}=0.
\end{eqnarray}

We note that these equations are different from those formulated for boundary layers in steady flow in \cite{RM3}. In the steady flow analysis of \cite{RM3}, it was assumed that there is no velocity boundary layer, i.e. $u$ is of order $1/W$ and $v$ is of order $1/W^2$, compared to order $1$ and $1/W$ above. Moreover, in a steady flow situation, the terms $\partial\S/\partial t$ and $(\v\cdot\nabla)\S$ both vanish on the boundary, and the term $\S/W$ becomes important. The parameter $E$ can be eliminated by a rescaling of $\S$, and we henceforth set it equal to one.

In an actual flow geometry, the boundary is curved. We shall show now that this does not change the boundary layer equations.
For simplicity, we shall stick to the two-dimensional case, where each component of the boundary is a closed curve. We shall use
local coordinates $q$ and $r$, where $q$ denotes arclength along the boundary, and $r$ is distance from the boundary.
In three dimensions, we would have to treat the boundary as a manifold, where coordinates can be defined only locally. However,
this would not fundamentally alter the analysis.

\def\i{{\bf i}}
\def\e{{\bf e}}

In three dimensions, let $x_i$ denote Cartesian coordinates, and let $p_i$ denote curvilinear but orthogonal coordinates.
We use some transformation rules given in Appendix I of \cite{Som}.
Let $g_k$ be defined by
\begin{equation}
g_k^2=\sum_i({\partial x_i\over\partial p_k})^2
\end{equation}
and let $\i_k, \e_k$ be the unit vectors in the Cartesian and $p$ coordinate systems, respectively.
Then we have
\begin{equation}
{\partial\e_k\over \partial p_l }={1\over g_k}{\partial g_l\over\partial p_k}\e_l-\delta_{kl}\sum_h{1\over g_h}{\partial g_k\over\partial p_h}\e_h,
\end{equation}
and
\begin{equation}
\nabla=\sum_k\e_k{1\over g_k}{\partial \over\partial p_k}.
\end{equation}

Hence for a scalar function $f$ we have the gradient
\begin{equation}
\nabla f=f_{x_1}\i_1+f_{x_2}\i_2+f_{x_3}\i_3=\sum\limits_{k=1}^3 {1\over g_k}{\partial f\over\partial p_k}\e_k,
\end{equation}
and the divergence of a vector ${\bf f}=f_1\e_1+f_2\e_2+f_3\e_3$ in the p system is given by
\begin{equation}
\nabla\cdot{\bf f}={1\over g_1g_2g_3}({\partial (f_1g_2g_3)\over\partial p_1}+{\partial (f_2g_1g_3)\over\partial p_2}+{\partial (f_3g_1g_2)\over\partial p_3}).
\end{equation}
In the two-dimensional case, this reduces to
\begin{equation}
\nabla\cdot{\bf f}={1\over g_1g_2}({\partial (f_1g_2)\over\partial p_1}+{\partial (f_2g_1)\over\partial p_2}).
\end{equation}

Convective terms in the equations transform as follows:
\begin{equation}
({\bf a}\cdot \nabla){\bf b}=\sum_{i,k}\e_k({a_i\over g_i}{\partial b_k\over\partial p_i}+{a_kb_i\over g_kg_i}{\partial g_k\over\partial p_i}-{a_ib_i\over g_ig_k}{\partial g_i\over\partial p_k}).
\end{equation}

Now we consider a solid boundary of our domain parameterized by the arc length $q$, $\x^*=\x^*(q)\in C^2$. In the local coordinates $(q,r)$ consider a point $\x=\x^*(q)+r\n$ Then we have
\begin{equation}
\left( \begin{array}{c} {\x^*_1}'(q) \\ {\x^*_2}'(q) \end{array} \right)=\left( \begin{array}{c} \n_2(q) \\ -\n_1(q)\end{array} \right), \left( \begin{array}{c} \n_1'(q) \\ \n_2'(q) \end{array} \right)=-\rho(q)\left( \begin{array}{c} \n_2(q) \\ -\n_1(q)\end{array} \right).
\end{equation}
Here $\rho(q)$ is the curvature which is assumed to be bounded.
When $r$ is small, ${\partial(\x_1,\x_2)\over \partial (q,r)}=1-\rho(q) r\neq 0$, hence $(q,r)$ are local coordinates.
For this coordinates transformation we have
\begin{eqnarray}
g_1&=&\sqrt{\x_1'(q)^2+\x_2'(q)^2}=1-\rho r,\nonumber\\
g_2&=&1.
\end{eqnarray}

The unit vectors in the curvilinear coordinate system are $\e_1=\tau, \e_2=\n$. With
$\v=u\e_1+v\e_2$, we obtain
\begin{equation}
\nabla\cdot\v={1\over 1-\rho r}{\partial u\over\partial q}+{\partial v\over\partial r}-{\rho\over 1-\rho r}v
\end{equation}
and
\begin{eqnarray}
(\v\cdot\nabla)\v&=&[{1\over 1-\rho r}u{\partial u\over\partial q}+v{\partial u\over\partial r}-{\rho uv\over 1-\rho r}]\e_1
\nonumber\\
&+&[{u\over 1-\rho r}{\partial v\over\partial r}+v{\partial v\over\partial r}+{\rho u^2\over 1-\rho r}]\e_2.
\end{eqnarray}

For the divergence of the stress tensor we find, with $\S=\sum_{ij} S_{ij}\e_i\e_j$,
\begin{eqnarray}
\nabla\cdot\S&=&[{1\over 1-\rho r}{\partial S_{11}\over\partial q}+{\partial S_{12}\over\partial r}-{2\rho\over 1-\rho r}S_{12}]\e_1\nonumber\\
&+& [{1\over 1-\rho r}{\partial S_{12}\over\partial q}+{\partial S_{22}\over\partial r}+{\rho\over 1-\rho r}S_{11}-{\rho\over 1-\rho r}S_{22}]\e_2.\end{eqnarray}

Under the scaling
\begin{eqnarray}
q'=q, r'=Wr, u'=u, v'=Wv,\nonumber\\
S_{11}'=S_{11}, S_{12}'=WS_{12},S_{22}'=W^2S_{22}, p'=p,
\end{eqnarray}
the divergence condition becomes
\begin{equation}
\nabla\cdot\v'={1\over 1-{r'\over W}\rho}{\partial\over\partial q'}u'+{\partial {v'\over W}\over\partial {r'\over W}}+({-\rho\over 1-{r'\over W}\rho}{v'\over W})=0.
\end{equation}
To the leading order as $W\rightarrow\infty$, we find
\begin{equation}
\nabla\cdot\v'={\partial u'\over\partial q'}+{\partial v'\over\partial r'}=0.
\end{equation}
That is, the divergence condition assumes the same form in the $(q',r')$ coordinates as in (\ref{scaling2}) above.
In a similar fashion, it can be seen that the momentum and constitutive equations also remain unchanged at leading order.
The study of the boundary layer problem therefore reduces to finding spatially periodic solutions (with respect to $x$) of
the equations (\ref{scaling1}) and (\ref{scaling2}) above.

\section{Lagrangian description}

We shall transform (\ref{scaling1}) and (\ref{scaling2}) into Lagrangian coordinates. Let $\xi_1$ and $\xi_2$ denote the coordinates of a fluid particle at $t=0$, and let
$x=x(\xi_1,\xi_2,t)$, $y=y(\xi_1,\xi_2,t)$ denote the position of the same particle at a later time.
That is, we have
\begin{eqnarray}
{\partial x(\xi_1,\xi_2,t)\over\partial t}&=&u(x(\xi_1,\xi_2,t),y(\xi_1,\xi_2,t),t),\nonumber\\
{\partial y(\xi_1,\xi_2,t)\over\partial t}&=&v(x(\xi_1,\xi_2,t),y(\xi_1,\xi_2,t),t),\\
x(\xi_1,\xi_2,0)&=&\xi_1,\nonumber\\
y(\xi_1,\xi_2,0)&=&\xi_2.\nonumber
\end{eqnarray}

We introduce the deformation gradient tensor $\F=\left( \begin{array}{cc} {\partial x\over\partial \xi_1} & {\partial x\over\partial \xi_2} \\
{\partial y\over\partial \xi_1} & {\partial y\over\partial \xi_2}\end{array} \right)$. It is easy to check that it satisfies ${\partial \F\over \partial t}+(\v\cdot\nabla)\F-(\nabla\v)\F=0,\F(0)=I$ and $\F$ is non-singular all the time.
Let $\C=\F^{-1}\S\F^{-T}$. It then follows from the constitutive law that ${\partial\C\over\partial t}+(\v\cdot\nabla)\C=0$.

The boundary layer system now becomes
\begin{eqnarray}\label{mainC}
{\partial u\over \partial t}+u{\partial u\over \partial x}+v{\partial u\over \partial y}={\partial S_{11}\over \partial x}+{\partial S_{12}\over \partial y}-P(x,t),\nonumber\\
{\partial\C\over\partial t}+(\v\cdot\nabla)\C=\0 ,\nonumber\\
{\partial u\over \partial x}+{\partial v\over \partial y}=0,\\ \nonumber
u=0, v=0, \mbox{at }  y=0,\\
u(x,y,0)=u_0(x,y), \nonumber\\
\C(x,y,0)=\S_0(x,y).\nonumber
\end{eqnarray}
Here we have introduced $P(x,t)$ for the known pressure gradient $\partial p/\partial x$.

We find that in Lagrangian coordinates, we have
${\partial\C\over\partial t}=0$, and hence\break $\C(\xi_1,\xi_2,t)=\C(\xi_1,\xi_2)$. We shall assume throughout that $\C$ is positive semidefinite and  there is positive constant $C_0$ such that $C_{22}\geq C_0>0$.

The following lemma shows how our
governing equations transform.

{\bf Lemma 1.}  In Lagrangian coordinates, from (\ref{mainC}) we have for $x$ the second order hyperbolic equation:
\begin{eqnarray}\label{mainL}
{\partial^2 x\over \partial t^2}&=&\nabla\cdot({\bf C}\nabla x)-P(x,t),\nonumber\\
x(\xi_1,\xi_2,0)&=&\xi_1,\nonumber\\
{\partial x\over\partial t}(\xi_1,\xi_2,0)&=&u_0(\xi_1,\xi_2),\\
x(\xi_1,0,t)&=&\xi_1. \nonumber
\end{eqnarray}

{\bf Proof:}
The horizontal acceleration is given by
\begin{equation}
{\partial u\over \partial t}+u{\partial u\over \partial x}+v{\partial u\over \partial y}={\partial^2 x(\xi_1,
\xi_2,t)\over \partial t^2}.
\end{equation}
Next we show that
\begin{equation}
{\partial S_{11}\over \partial x}+{\partial S_{12}\over \partial y}={\partial \over \partial \xi_i}(C_{ji}(\xi_1,\xi_2){\partial x\over\partial \xi_j})=\nabla\cdot({\bf C}\nabla x).\end{equation}

We write $\divg(\F\C\F^T)$ in components, with the Einstein summation convention:
\begin{equation}
{\partial\over\partial x_j}(F_{ik}C_{kl}F_{jl})=F_{jl}{\partial\over\partial x_j}(C_{kl}F_{ik})+
C_{kl}F_{ik}{\partial F_{jl}\over\partial x_j}.
\end{equation}
For the first term, we note that by the chain rule, we have
\begin{equation}
F_{jl}{\partial\over\partial x_j}={\partial x_j\over\partial\xi_l}{\partial\over\partial x_j}={\partial\over\partial\xi_l}.
\end{equation}
For the second term, we write
\begin{equation}
{\partial F_{jl}\over\partial x_j}={\partial \xi_k\over\partial x_j}{\partial F_{jl}\over\partial \xi_k}=
{\partial \xi_k\over\partial x_j}{\partial^2x_j\over\partial\xi_k\partial x_l}={\rm tr}\,(\F^{-1}{\partial\F\over\partial\xi_l}).
\end{equation}

Since ${\rm det}\,\F=1$ by the incompressibility condition, we find
\begin{equation}
{\partial\over\partial \xi_l}({\rm det}\,\F)=({\rm det}\,\F){\rm tr}\,(\F^{-1}{\partial\F\over\partial\xi_l})=0.
\end{equation}

Once we have found $x$, we can solve for $y$ from the incompressibility condition.

{\bf Lemma 2.} From the incompressibility condition, we find the following equation for $y$ under the Lagrangian description:
\begin{eqnarray}\label{main_continuity}
{\partial x\over\partial\xi_1}{\partial y\over\partial \xi_2}-{\partial x\over\partial\xi_2}{\partial y\over\partial \xi_1}=1,\nonumber \\
y(\xi_1,0,t)=0.
\end{eqnarray}

This equation is a well-posed transport equation for $y$ with $t$ as a parameter. The details are omitted.

After we get $x(\xi_1,\xi_2,t)$ and $y(\xi_1,\xi_2,t)$, we can recover $u=x_t, v=y_t$, hence the deformation tensor $\F=\left( \begin{array}{cc} {\partial x\over\partial \xi_1} & {\partial x\over\partial \xi_2} \\
{\partial y\over\partial \xi_1} & {\partial y\over\partial \xi_2}\end{array} \right)$ and extra stress $\S=\F\C\F^T$.

\section{Well-posedness result}

We now focus on the solution of (\ref{mainL}). This is simply a semilinear wave equation. We are seeking solutions which are periodic with respect
to $\xi_1$ with a given period, but the behavior for $\xi_2\to\infty$ warrants some discussion. We want to match to a given outer solution. Therefore,
we assume that the initial condition $u_0$ has a limit at infinity:
\begin{equation}\label{u0uinfinity}
\lim_{\xi_2\to\infty} u_0(\xi_1,\xi_2)=u^\infty(\xi_1).
\end{equation}

We also assume that $C_{11}$ has a limit,
\begin{equation}
\lim_{\xi_2\to\infty} C_{11}(\xi_1,\xi_2)=C_{11}^\infty(\xi_1),
\end{equation}
and that $C_{12}$ and $C_{22}$ and their derivatives are uniformly bounded.

The behavior at infinity is now governed by the following limit problem:
\begin{eqnarray}\label{outflow}
{\partial^2x^\infty(\xi_1,t)\over\partial t^2}&=&{\partial\over\partial\xi_1}(C_{11}^\infty(\xi_1){\partial x^\infty\over\partial\xi_1})-P(x^\infty,t),\nonumber\\
x^\infty(\xi_1,0)&=&\xi_1,\nonumber\\
{\partial x^\infty\over\partial t}(\xi_1,0)&=&u^\infty(\xi_1).
\end{eqnarray}
Note that this is precisely the $x$ component of the momentum equation for the outer flow at the wall, transformed to Lagrangian coordinates.

Combine (\ref{mainL}) and (\ref{outflow}), and denote $x-x^\infty$ by $X$. We have
\begin{eqnarray}\label{mainL1}
X_{tt}&=&\nabla\cdot({\bf C}\nabla X)+ \Psi(X,\xi_1,\xi_2,t),\nonumber\\
X(\xi_1,\xi_2,0)&=&0,\nonumber\\
{\partial X\over\partial t}(\xi_1,\xi_2,0)&=&{\bf f}(\xi_1,\xi_2),\\
X(\xi_1,0,t)&=&{\bf g}(\xi_1,t) \nonumber
\end{eqnarray}
with
\begin{eqnarray}\label{mainL1data}
{\Psi}&=&\nabla\cdot({\bf C}\nabla x^{\infty})-{\partial \over \partial \xi_1}(C_{11}^\infty(\xi_1){\partial x^\infty\over\partial \xi_1})-P(X+x^\infty,t)+P(x^\infty,t),\nonumber\\
{\bf f}(\xi_1,\xi_2)&=&u_0(\xi_1,\xi_2)-u^{\infty}(\xi_1),\\
{\bf g}(\xi_1,t)&=&\xi_1-x^\infty(\xi_1,t).\nonumber
\end{eqnarray}
Now we want to make the boundary and initial conditions homogeneous.
Pick a smooth function
$\chi=\left\{ \begin{array}{rcl}
1, & 0\leq\xi_2\leq 1,\\
0, & 2\leq \xi_2,
\end{array}\right. $  defined on $[0,\infty)$ such that $0\leq \chi\leq 1$.

Let
\begin{eqnarray}\label{X2x}
Y&=&{\bf g}(\xi_1,t)\cdot\chi(\xi_2)+t[u_0(\xi_1,\xi_2)-u^\infty(\xi_1)\cdot(1-\chi(\xi_2))],\nonumber\\
\bar{x}&=&X(\xi_1,\xi_2,t)-Y.
\end{eqnarray}

We have for $\bar{x}$
\begin{eqnarray}\label{mainL11}
\bar{x}_{tt}&=&\nabla\cdot({\bf C}\nabla \bar{x})+{\Phi}(\bar{x},\xi_1,\xi_2,t),\nonumber\\
\bar{x}(\xi_1,\xi_2,0)&=&0,\nonumber\\
{\partial \bar{x}\over\partial t}(\xi_1,\xi_2,0)&=&0,\\
\bar{x}(\xi_1,0,t)&=&0, \nonumber
\end{eqnarray}
where
\begin{equation}\label{Phi}
{\Phi}(\bar{x},\xi_1,\xi_2,t)= \Psi(\bar{x}+Y,\xi_1,\xi_2,t)+ \nabla\cdot({\bf C}\nabla Y)-{Y}_{tt}.
\end{equation}

We have an initial-boundary value problem for a semilinear wave equation of $\bar{x}$. We shall state our existence result for $\bar{x}$ in $L^2$ based spaces.  Unlike the work of Lasiecka et al. \cite{LLT}, our problem is posed on an unbounded domain and does not satisfy the uniformly hyperbolic condition, namely we do not require $\C>0$ strictly throughout the domain. To study this IBVP, we shall transform the equation into an equivalent first order system.

Since $\C$ is real positive semi-definite, there is a unique real positive semi-definite $\bf A$ with $a_{22}\geq c_0>0$, denoted by ${\bf A}=\left( \begin{array}{cc} a_{11} & a_{12} \\
a_{12} & a_{22} \end{array} \right)$,  such that ${\bf A}^2=\C$. Let $U=\bar{x}_t, (V,W)^T={\bf A}\nabla \bar{x}, \vec{V}=(U,V,W)^T$. >From (\ref{mainL11}) we have an integro-differential system
\begin{eqnarray}\label{mainL21}
\vec{V}_t
&=&{\bf A}_1 \vec{V}_{\xi_1}
+{\bf A}_2\vec{V}_{\xi_2}+{\bf \B}\vec{V}+\vec{\Phi},
\end{eqnarray}
with initial and boundary condition
\begin{eqnarray}\label{mainL21_IBC}
\vec{V}(\xi_1,\xi_2,0)&=& {\vec{0}}, \nonumber\\
U(\xi_1,0,t)&=&\bar{x}_t(\xi_1,0,t)=0.
\end{eqnarray}

Here  ${\bf A}_1=\left( \begin{array}{ccc}
0& a_{11} & a_{12}\\
a_{11} & 0 & 0 \\
a_{12} & 0 & 0 \end{array} \right),
{\bf A}_2= \left( \begin{array}{ccc}
0& a_{12} & a_{22}\\
a_{12} & 0 & 0 \\
a_{22} & 0 & 0 \end{array} \right),
\vec{\Phi}=\left( \begin{array}{ccc}
\Phi\\
0  \\
0 \end{array} \right),
$

$ {\B}= \left( \begin{array}{ccc}
0& {\partial a_{11}\over {\partial \xi_1}}+{\partial a_{12}\over \partial{\xi_2}} &  {\partial a_{12}\over {\partial \xi_1}}+{\partial a_{22}\over \partial{\xi_2}}\\
0 & 0 & 0 \\
0 & 0 & 0 \end{array} \right), \Phi=\Phi([\int_0^tU(\xi_1,\xi_2,s) ds], \xi_1,\xi_2,t).$

This is a first order linear symmetric hyperbolic system when $\Phi$ is given. The dependence of $\Phi$ on $\bar{x}$ can then be handled via a standard fixed point argument.
In the following, we shall show that on the half plane with given periodic forcing term $\Phi$, the system (\ref{mainL21})-(\ref{mainL21_IBC}) has a unique local solution which is periodic.

The theory for hyperbolic system was developed by Friedrichs \cite{KOF}, Kreiss \cite{KHO}, Lax and Phillips \cite{LP1}, and others. The general idea is to do energy estimate on the equations. Based on proper a priori estimates, one can either define a weak solution then improve its regularity \cite{KOF0}, or use different schemes to approximate the real solutions \cite{KL}. For general hyperbolic system with characteristic boundary, the full regularity may be lost. People developed anisotropic weighted space of Sobolev type to meet this purpose, see Secchi \cite{SP}.

Our problem presents a first order symmetric hyperbolic system with characteristic boundary. We show that nevertheless we are still able to get the full regularity in Sobolev spaces because of its special structure.
\def\enn{{\rm I\!N}}

{\bf Theorem 1.}
Consider the boundary layer system (\ref{mainL21})-(\ref{mainL21_IBC}) on the domain
$$-\infty<\xi_1<\infty,  0\leq \xi_2<\infty, \mbox{ for   }  t\geq 0.$$

Let $m\ge 1$ be an integer. We assume that $\A_1,\A_2, \B$ and their derivatives up to order $m$ are bounded for $(\xi_1,\xi_2,t)\in (-\infty,\infty)\times[0,\infty)\times [0,T]$. Moreover they are all periodic in $\xi_1$ with period $1$. Let $\Omega =[0,1]\times [0,\infty) $, $ \Q=\Omega\times [0,T]$, and denote by $H_p^{m}(\Q), m\in \enn$ the space of all periodic (in $\xi_1$) functions which have $H^m$ regularity.
If $\Phi \in H_p^{m}(\Q)$, then there exists some $T'\in(0,T]$ such that a unique solution of (\ref{mainL21})-(\ref{mainL21_IBC})  exists and satisfies $\vec{V}\in H_p^{m}(\Q)$.

The idea of the proof is to develop a priori estimates for solutions which are assumed smooth, a density argument can be used to extend the result for general Sobolev data in $H^m$.
In the estimates, the inhomogeneous forcing term $\Phi$ always contributes terms of the form $(D^\alpha\Phi, D^\beta U)\leq C (||D^\alpha\Phi||^2+||D^\beta U||^2)$. Since these terms
are easily dealt with, we present the estimates without the forcing term $\Phi$.

{\bf Proof:}
We complete the energy estimate in an elementary way. The $L^2$ estimate is done by virtue of the special structure of coefficient matrices. The evaluation of higher order Sobolev norms gets into trouble due to the loss of normal derivative at the boundary. For that reason we do the estimate separately. We first work on the $L^2$ norm. Then we use the integral for the estimate of all derivative except the normal direction. For the normal derivative, we take advantage of the boundary condition and of a constraint on the solution that is preserved by the equations.

{\bf 1.  $L^2$ estimate}

Consider (\ref{mainL21})-(\ref{mainL21_IBC}).
Note that boundary conditions are imposed only on $U$, and not on $V$ and $W$.

We multiply (\ref{mainL21}) with $\vec{V}$ and integrate it with respect to $\xi\in\Omega$
\begin{eqnarray}\label{estimate2}
\int_\Omega\vec{V}_t\cdot \vec{V}d\xi =\int_{\Omega}{\bf A}_1\vec{V}_{\xi_1}\cdot \vec{V}d\xi+\int_{\Omega}{\bf A}_2\vec{V}_{\xi_2}\cdot \vec{V}d\xi+\int_{\Omega}\B\vec{V}\cdot \vec{V}d\xi.
\end{eqnarray}

Integrating by parts, we have

\begin{eqnarray}\label{estimate3}
{1\over 2}{d\over dt}\int_\Omega|\vec{V}|^2d\xi
&=&{1\over 2}\int_{0}^{\infty}{\bf A}_1\vec{V}\cdot \vec{V}|_{\xi_1=0}^1d\xi_2-{1\over 2} \int_{\Omega}{\partial {\bf A}_1\over\partial \xi_1}\vec{V}\cdot \vec{V}d\xi\nonumber\\
&  &+
{1\over 2}\int_{0}^{1}{\bf A}_2\vec{V}\cdot \vec{V}|_{\xi_2=0}^{\infty}d\xi_1-{1\over 2} \int_{\Omega}{\partial {\bf A}_2\over\partial \xi_2}\vec{V}\cdot \vec{V}d\xi\nonumber\\
& & +{1\over 2}\int_{\Omega}(\B+\B^*)\vec{V}\cdot \vec{V}d\xi\nonumber\\
&=& {1\over 2}\int_{\Omega}[\B+\B^*-{\partial {\bf A}_1\over\partial \xi_1}-{\partial {\bf A}_2\over\partial \xi_2}]\vec{V}\cdot \vec{V}d\xi.
\end{eqnarray}

because $\vec{V}(0,\xi_2,t)=\vec{V}(1,\xi_2,t)$ and $\int_{0}^{1}{\bf A}_2\vec{V}(\xi_1,0,t)\cdot \vec{V}(\xi_1,0,t)d\xi_1=0$.

{\bf 2.  $\partial_t$ estimate}

The energy estimate for $ \vec{V}_t$ can be done similarly to the $L^2$ case. Differentiating (\ref{mainL21}) w.r.t $t$ we have
\begin{eqnarray}\label{estimate_t1}
(\vec{V}_{t})_t
={\bf A}_1 (\vec{V}_t)_{\xi_1} +{\bf A}_2(\vec{V}_{t})_{\xi_2}+\B\vec{V}_{t}.
\end{eqnarray}

Multiply and integrate to have
\begin{eqnarray}\label{estimate_t2}
\int (\vec{V}_t)_t \cdot \vec{V}_{t}d\xi
=\int {\bf A}_1 (\vec{V}_t)_{\xi_1}\cdot \vec{V}_td\xi\nonumber\\
+\int {\bf A}_2(\vec{V}_t)_{\xi_2}\cdot \vec{V}_td\xi+\int \B\vec{V}_{t}\cdot \vec{V}_{t}d\xi.
\end{eqnarray}

Simplify to have
\begin{eqnarray}\label{estimate_t3}
{1\over 2}{d\over dt}\int_\Omega|\vec{V}_t|^2d\xi
&=& -{1\over 2}\int_{\Omega}{\partial {\bf A}_1\over\partial \xi_1}\vec{V}_{t}\cdot \vec{V}_t d\xi + {1\over 2}\int_{0}^{1}{\bf A}_2\vec{V}_t\cdot \vec{V}_t|_{\xi_2=0}^{\infty}d\xi_1\nonumber\\
&  &- {1\over 2}\int_{\Omega}{\partial {\bf A}_2\over\partial \xi_2}\vec{V}_t\cdot \vec{V}_td\xi+{1\over 2}\int_{\Omega}(\B+\B^*)\vec{V}_t\cdot \vec{V}_td\xi\nonumber\\
&=& {1\over 2}\int_{\Omega}[\B+\B^*-{\partial {\bf A}_1\over\partial \xi_1}-{\partial {\bf A}_2\over\partial \xi_2}]\vec{V}_t\cdot \vec{V}_td\xi.
\end{eqnarray}
Notice that ${\bf A}_2\vec{V}_t\cdot\vec{V}_t=0$ at the boundary.

{\bf 3.  $\partial_{\xi_1}$ estimate}

Differentiate (\ref{mainL21}) w.r.t $\xi_1$ we have
\begin{eqnarray}\label{estimate_xi11}
(\vec{V}_{\xi_1})_t
={\bf A}_1 (\vec{V}_{\xi_1})_{\xi_1}+{\partial {\bf A}_1\over \partial {\xi_1}} \vec{V}_{\xi_1}\nonumber\\
+{\bf A}_2(\vec{V}_{\xi_1})_{\xi_2}+{\partial {\bf A}_2\over \partial {\xi_1}}\vec{V}_{\xi_2}+\B \vec{V}_{\xi_1}+{\partial \B\over\partial \xi_1 }\vec{V}.
\end{eqnarray}
Multiply and integrate to have

\begin{eqnarray}\label{estimate_xi12}
\int (\vec{V}_{\xi_1})_t \cdot (\vec{V}_{\xi_1})d\xi
&=& \int {\bf A}_1 (\vec{V}_{\xi_1})_{\xi_1}\cdot \vec{V}_{\xi_1}d\xi +\int {\partial {\bf A}_1\over \partial {\xi_1}} \vec{V}_{\xi_1}\cdot \vec{V}_{\xi_1}d\xi\nonumber\\
& & +\int {\bf A}_2(\vec{V}_{\xi_1})_{\xi_2}\cdot \vec{V}_{\xi_1}d\xi+\int {\partial {\bf A}_2\over \partial {\xi_1}}\vec{V}_{\xi_2}\cdot \vec{V}_{\xi_1}d\xi\nonumber\\
& & +\int \B \vec{V}_{\xi_1}\cdot \vec{V}_{\xi_1}d\xi+\int {\partial \B\over\partial \xi_1 }\vec{V}\cdot \vec{V}_{\xi_1}d\xi.
\end{eqnarray}

Simplify to have
\begin{eqnarray}\label{estimate_xi13}
{1\over 2}{d\over dt}\int_\Omega|\vec{V}_{\xi_1}|^2d\xi
&=& -{1\over 2}\int_{\Omega}{\partial {\bf A}_1\over\partial \xi_1}\vec{V}_{\xi_1}\cdot \vec{V}_{\xi_1}d\xi+\int {\partial {\bf A}_1\over \partial {\xi_1}} \vec{V}_{\xi_1}\cdot \vec{V}_{\xi_1}d\xi\nonumber\\
&  &+ {1\over 2}\int_{0}^{1}{\bf A}_2\vec{V}_{\xi_1}\cdot \vec{V}_{\xi_1}|_{\xi_2=0}^{\infty}d\xi_1- {1\over 2}\int_{\Omega}{\partial {\bf A}_2\over\partial \xi_2}\vec{V}_{\xi_1}\cdot \vec{V}_{\xi_1}d\xi\nonumber\\
&  &  +\int {1\over 2}[\B+\B^*] \vec{V}_{\xi_1}\cdot \vec{V}_{\xi_1}d\xi\nonumber\\
& &+{1\over 2}\int {\partial (\B+\B^*)\over\partial \xi_1 }\vec{V}\cdot \vec{V}_{\xi_1}d\xi+ \int {\partial {\bf A}_2\over \partial {\xi_1}}\vec{V}_{\xi_2}\cdot \vec{V}_{\xi_1}d\xi\nonumber\\
&=& {1\over 2}\int_{\Omega}[{\partial {\bf A}_1\over\partial \xi_1}-{\partial {\bf A}_2\over\partial \xi_2}+\B+\B^* ]\vec{V}_{\xi_1}\cdot \vec{V}_{\xi_1}d\xi\nonumber\\
& & +{1\over 2}\int {\partial (\B+\B^*)\over\partial \xi_1 }\vec{V}\cdot \vec{V}_{\xi_1}d\xi+ \int {\partial {\bf A}_2\over \partial {\xi_1}}\vec{V}_{\xi_2}\cdot \vec{V}_{\xi_1}d\xi.
\end{eqnarray}

To estimate $\vec{V}_{\xi_1}$, we need knowledge of $\vec{V}_{\xi_2}$ because of the term
$$ \int {\partial {\bf A}_2\over \partial {\xi_1}}\vec{V}_{\xi_2}\cdot \vec{V}_{\xi_1}d\xi.$$
This turns out to be no problem when we finish the {$\partial_{\xi_2}$ estimate} below.

{\bf 4.  $\partial_{\xi_2}$ estimate}

The crucial estimate is for the term $\vec{V}_{\xi_2}$. Because of the loss of normal derivative, we can not proceed as before. Instead, we go back to original equation and solve for $\vec{V}_{\xi_2}$ as a function of $\vec{V}, \vec{V}_t, \vec{V}_{\xi_1}$.

>From (\ref{mainL21}) we get
\begin{eqnarray}\label{original11}
U_t=a_{12}V_{\xi_2}+a_{22}W_{\xi_2} +a_{11}V_{\xi_1}+a_{12}W_{\xi_1}\\ \nonumber
+ [{\partial a_{11}\over {\partial \xi_1}}+{\partial a_{12}\over \partial{\xi_2}}]V+[ {\partial a_{12}\over {\partial \xi_1}}+{\partial a_{22}\over \partial{\xi_2}}]W,
\end{eqnarray}
\begin{eqnarray}\label{original12}
V_t=a_{11}U_{\xi_1}+a_{12}U_{\xi_2},
\end{eqnarray}
\begin{eqnarray}\label{original13}
W_t=a_{12}U_{\xi_1}+a_{22}U_{\xi_2}.
\end{eqnarray}

By (\ref{original13}) we have $$U_{\xi_2}={1\over a_{22}}(W_t-a_{12}U_{\xi_1}),$$ which leads to
$$||U_{\xi_2}||_0\leq C (||W_t||_0+||U_{\xi_1}||_0).$$

To solve for $V_{\xi_2}$ and $W_{\xi_2}$, we differentiate (\ref{original12}) with respect to $\xi_1$, $\xi_2$, respectively, to have

$$V_{\xi_1 t}=a_{11}U_{\xi_1\xi_1}+a_{12}U_{\xi_1\xi_2}+{\partial a_{11}\over\partial \xi_1 }U_{\xi_1}+{\partial a_{12}\over\partial \xi_1 }U_{\xi_2},$$

$$V_{\xi_2 t}=a_{11}U_{\xi_1\xi_2}+a_{12}U_{\xi_2\xi_2}+{\partial a_{11}\over\partial \xi_2 }U_{\xi_1}+{\partial a_{12}\over\partial \xi_2 }U_{\xi_2}.$$

Similarly, for (\ref{original13}) we have

$$W_{\xi_1 t}=a_{12}U_{\xi_1\xi_1}+a_{22}U_{\xi_1\xi_2}+{\partial a_{12}\over\partial \xi_1 }U_{\xi_1}+{\partial a_{22}\over\partial \xi_1 }U_{\xi_2},$$

$$W_{\xi_2 t}=a_{12}U_{\xi_1\xi_2}+a_{22}U_{\xi_2\xi_2}+{\partial a_{12}\over\partial \xi_2 }U_{\xi_1}+{\partial a_{22}\over\partial \xi_2 }U_{\xi_2}.$$

Multiplying with $a_{12}, a_{22}$ and doing simple calculations, we obtain

$$(a_{22}V_{\xi_2}-a_{12}W_{\xi_2}+a_{12}V_{\xi_1}-a_{11}W_{\xi_1})_t$$

$$=[a_{12}{\partial a_{11}\over \partial \xi_1}+a_{22}{\partial a_{11}\over\partial \xi_2}-a_{11}{\partial a_{12}\over\partial \xi_1}-a_{12}{\partial a_{12}\over\partial \xi_2}]U_{\xi_1}$$

$$+[a_{12}{\partial a_{12}\over \partial \xi_1}+a_{22}{\partial a_{12}\over\partial \xi_2}-a_{11}{\partial a_{22}\over\partial \xi_1}-a_{12}{\partial a_{22}\over\partial \xi_2}]U_{\xi_2}$$

$$=[a_{12}{\partial a_{11}\over \partial \xi_1}+a_{22}{\partial a_{11}\over\partial \xi_2}-a_{11}{\partial a_{12}\over\partial \xi_1}-a_{12}{\partial a_{12}\over\partial \xi_2}]U_{\xi_1}$$

$$+{1\over a_{22}}[a_{12}{\partial a_{12}\over \partial \xi_1}+a_{22}{\partial a_{12}\over\partial \xi_2}-a_{11}{\partial a_{22}\over\partial \xi_1}-a_{12}{\partial a_{22}\over\partial \xi_2}](W_t-a_{12}U_{\xi_1})$$

$$=\F(U_{\xi_1}, W_t, \xi_1,\xi_2).$$

Since $\vec{V}(\xi_1,\xi_2,0)=0$,  we have
\begin{eqnarray}\label{original14}
a_{22}V_{\xi_2}-a_{12}W_{\xi_2}=\int_0^t \F(U_{\xi_1}, W_t, \xi_1,\xi_2)ds-a_{12}V_{\xi_1}+a_{11}W_{\xi_1}.
\end{eqnarray}

And (\ref{original11}) yields
\begin{eqnarray}\label{original15}
a_{12}V_{\xi_2}+a_{22}W_{\xi_2}=U_t-a_{11}V_{\xi_1}-a_{12}W_{\xi_1}\\ \nonumber
-[{\partial a_{11}\over {\partial \xi_1}}+{\partial a_{12}\over \partial{\xi_2}}]V+[ {\partial a_{12}\over {\partial \xi_1}}+{\partial a_{22}\over \partial{\xi_2}}]W.
\end{eqnarray}

Together we can now solve for $V_{\xi_2}$ and $W_{\xi_2}$ in terms of $U_t$, $U_{\xi_1}$, $V_{\xi_1}$, $W_{\xi_1}$, $V$, $W$, and we know $\vec{V}_{\xi_2}$ in term of other derivatives. We then plug  $\vec{V}_{\xi_2}$ back in (\ref{estimate_xi13}) to complete the estimate for $\vec{V}_{\xi_1}$.  Finally, the estimate for $\vec{V}_{\xi_2}$ is done by using the triangle inequality.  For the estimate of the integral term  $\int_0^t \F(U_{\xi_1}, W_t, \xi_1,\xi_2)ds$, we use Minkowski's inequality. In this way we get the $H_1$ regularity of the solution. Taking higher order derivatives of (\ref{mainL21}), the $H_m$ energy estimate can be done in same fashion.

Based on the energy estimate, the uniqueness and regularity are easy to see. The existence proof can be obtained along the lines of \cite{KL}. Our system corresponds to Case 2 in the discussion
given in that reference. Following \cite{KL}, to show the existence, we would first diagonalize the system so that $\A_2$ transforms to
$$\left( \begin{array}{ccc}
-\lambda& 0 & 0\\
0 & 0 & 0 \\
0 & 0 & +\lambda \end{array} \right),$$
where $\lambda=\sqrt{a_{12}^2+a_{22}^2}\geq c_0>0$. Then we consider the perturbed problem with
$$ \A_{2\sigma}=\left( \begin{array}{ccc}
-\lambda+\sigma& 0 & 0\\
0 & \sigma & 0 \\
0 & 0 & +\lambda+\sigma \end{array} \right),$$
where $0<\sigma<\lambda$. This perturbed problem has same number of negative eigenvalue as the original problem, but a noncharacteristic boundary. Standard methods can be used to show that there exists a solution depending on $\sigma$, which, for instance, is constructed by finite difference approximating scheme in \cite{KL}. The estimates are independent of $\sigma$ and a limit argument gives the existence of original problem (\ref{mainL21}).  We refer to \cite{KL}, page 297-298 for a detailed discussion.

\end{document}